\documentclass[11pt,bezier]{article}
\usepackage{amsmath}
\usepackage{amsfonts,amsthm,amssymb}
\usepackage{amsfonts}
\usepackage{graphics}
\usepackage{amssymb}

\newtheorem{lem}{Lemma}[section]
\newtheorem{thm}[lem]{Theorem}

\newtheorem{conj}{Conjecture}

\theoremstyle{definition}

\begin{document}
\title{Connectivity keeping caterpillars and spiders in bipartite graphs with connectivity at most three\footnote{The research is supported by National Natural Science Foundation of China (11861066).}}
\author{Qing Yang, Yingzhi Tian\footnote{Corresponding author. E-mail: tianyzhxj@163.com (Y. Tian).} \\
{\small College of Mathematics and System Sciences, Xinjiang
University, Urumqi, Xinjiang, 830046, PR China}}

\date{}

\maketitle

\noindent{\bf Abstract } A conjecture of Luo, Tian and Wu (2022) says that for every positive integer $k$ and every finite tree $T$ with bipartition $X$ and $Y$ (denote $t = \max\{|X|,|Y |\})$, every $k$-connected bipartite graph $G$ with $\delta(G) \geq k + t$ contains a subtree $T' \cong T$ such that $\kappa(G-V (T')) \geq k$. In this paper, we confirm this conjecture for caterpillars when $k=3$ and spiders when $k\leq 3$.

\noindent{\bf Keywords:} Connectivity, Bipartite graphs, Caterpillars, Spiders

\section{Introduction}

Throughout this paper, graph always means a finite, undirected graph without multiple edges and loops. Let $G$ be a graph with vertex set $V(G)$ and edge set $E(G)$. The \emph{minimum degree} and \emph{order} of $G$ are denoted by $\delta(G)$ and $|V(G)|$, respectively. The \emph{floor} of a real number $x$, denoted by $\lfloor x\rfloor$, is the greatest integer not larger than $x$; the \emph{ceiling} of a real number $x$, denoted by $\lceil x \rceil$, is the least integer greater than or equal to $x$. For graph-theoretical terminologies and notation not defined here, we follow \cite{Bondy}.

The $connectivity$ of $G$, denoted by $\kappa(G)$, is the minimum size of a vertex set $S$ such that $G -S$ is disconnected or has only one vertex. A graph $G$ is \emph{$k$-connected} if $\kappa(G)\geq k$.

In 1972, Chartrand, Kaugars and Lick \cite{Chartrand} proved that every $k$-connected graph $G$ with $\delta(G) \geq \lfloor \frac{3k}{2}\rfloor $ has a vertex $v$ with $\kappa(G -\{v\})\geq k$. Fujita and Kawarabayashi \cite{Fujita} showed every $k$-connected graph $G$ with $\delta(G) \geq \lfloor \frac{3k}{2}\rfloor +2$ contains an edge $uv$ such that $G -\{u,v\}$ remains $k$-connected. Meanwhile, they conjectured that every $k$-connected graph $G$ with $\delta(G) \geq \lfloor \frac{3k}{2}\rfloor+ f_{k}(m)-1$ contains a connected subgraph $G_1$ of order $m$ such that $\kappa(G-V (G_1))\geq k$, where $f_{k}(m)$ is non-negative. Mader confirmed the conjecture and proved the following theorem in \cite{Mader1}.

\begin{thm} (Mader \cite{Mader1}) Every $k$-connected graph $G $ with $\delta(G) \geq \lfloor \frac{3k}{2}\rfloor + m- 1$ contains a path $P$ of order $m$ such that $G -V (P)$ remains $k$-connected.
 \end{thm}

In the same paper, Mader stated the following conjecture for trees.

\begin{conj}(Mader \cite{Mader1})
For every tree $T$ of order $m$, every $k$-connected graph $G$ with $\delta(G) \geq \lfloor \frac{3k}{2}\rfloor+ m-1$ contains a subtree $T'\cong T$ such that $G-V (T')$ is still $k$-connected.
\end{conj}

Concerning to this conjecture, Mader \cite{Mader2} also proved that, for any tree $T$ with order $m$, every $k$-connected graph $G$ with $\delta(G) \geq 2(k-1+m)^2+m-1$ contains a subtree $T'\cong T$ such that $G-V (T')$ is still $k$-connected. For $k=1$, the result in \cite{Diwan} showed that Mader's conjecture is true. For $k=2$, Mader's conjecture was verified when $T$ is listed as follows: stars, double-stars, path-stars and path-double-stars \cite{Tian1,Tian2}; trees with small internal vertices, small girth and quasi-monotone caterpillars \cite{Hasunuma1,Hasunuma2}; trees with diameter at most 4 \cite{Lu}; caterpillars and spiders \cite{Hong1}. Recently, Hong and Liu \cite{Hong} proved that Mader's conjecture holds when $k\leq3$.

For  bipartite graphs, Luo, Tian and Wu \cite{Luo} proved the following theorem, which relaxes the minimum degree bound in Theorem 1.1.

\begin{thm} (Luo, Tian and Wu \cite{Luo})
Every $k$-connected bipartite graph $G$ with $\delta(G) \geq k + m$ contains a path $P$ of order $m$ such that $G-V (P)$ remains $k$-connected.
 \end{thm}

Analogue to Conjecture 1, for bipartite graphs, the authors in \cite{Luo} proposed the following conjecture.

 \begin{conj} (Luo, Tian and Wu)
 For every tree $T$ with bipartition $X$ and $Y$ (denote $t = \max\{|X|, |Y |\}$), every $k$-connected bipartite graph $G$ with $\delta(G) \geq k + t$ contains a subtree $T'\cong T$ such that $\kappa(G -V (T')) \geq k$.

\end{conj}

Since $K_{k+t-1,k+t-1}$ contains no $T'\cong T$ such that $\kappa(G - V (T')) \geq k$, the bound $\delta(G) \geq k + t$ would be best possible when the conjecture were true. Zhang \cite{Zhang} confirmed Conjecture 2 for caterpillars when $k\leq 2$.

In the next section, preliminaries will be given. We will confirm Conjecture 2 for caterpillars when $k=3$ in Section 3 and verify Conjecture 2 for spiders when $k\leq 3$ in Section 4.

\section{Preliminaries}

The set of neighbours of a vertex $v$ in $G$ is denoted by $N_{G}(v)$ and $d_{G}(v)=|N_{G}(v)|$. For any subset $S\subseteq V(G)$, $G[S]$ is the subgraph induced by $S$, and $N_{G}(S)$ is the neighborhood of $S$ in $G-S$. A \emph{$(u, v)$-path} $P$ is a path with ends $u$ and $v$. Then we write $End(P) = \{u, v\}$, say $P$ is a \emph{$u$-path} or a \emph{$v$-path}. For any $u',v'\in V(P)$, $u'Pv'$ is the subpath of $P$ between $u'$ and $v'$. For two subsets $V_{1}$ and $V_{2}$ of $V (G)$, the $(V_{1}, V_{2})$-path is a path with one end in $V_{1}$, the other end in $V_{2}$, but internal vertices not in $V_{1}\cup V_{2}$. We use $(v, V_{2} )$-path for $(\{v\}, V_{2})$-path.

\begin{lem}
Let $G$ be a bipartite graph.  ($i$) If $\delta(G)\geq \lceil \frac{m}{2}\rceil$, then $G$ contains a path of order at least $m$. ($ii$) If $\delta(G)\geq \lfloor \frac{m}{2}\rfloor$, then $G$ contains a path of order at least $m-1$. ($iii$) If $\delta(G)\geq m$, then $G$ contains a path of order at least $2m$.

 \end{lem}\noindent{\bf Proof.} We will construct a path in $G$ with order $m$ when $\delta(G)\geq \lceil \frac{m}{2}\rceil$. The result is clear for $m=1$ or $m=2$. We assume $m\geq 3$ and $P_{0}=v_{0}v_{1}$ is a path of order 2 in $G$. Since $|N_{G}(v_{1})\backslash V(P_{0})|\geq \lceil \frac{m}{2}\rceil-1\geq 1$, we can extend $P_{0}$ to a longer path $P_{1}=v_{0}v_{1}v_{2}$ in $G$, where $v_{2}\in N_{G}(v_{1})\backslash V(P_{0})$. Note that $|N_{G}(v_{2})\backslash V(P_{1})|\geq \lceil \frac{m}{2}\rceil-1\geq 1$. We can extend $P_{1}$ to a longer path $P_{2}=v_{0}v_{1}v_{2}v_{3}$ in $G$, where $v_{3}\in N_{G}(v_{2})\backslash V(P_{1})$. If $N_{G}(v_{3})\backslash V(P_{2})= \emptyset$, then $\lceil \frac{m}{2}\rceil\leq 2$ and $ m\leq |V(P_{2})|$. If $N_{G}(v_{3})\backslash V(P_{2})\neq \emptyset$, then we can extend $P_{2}$ to a longer path $P_{3}=v_{0}\cdots v_{4}$ in $G$, where $v_{4}\in N_{G}(v_{3})\backslash V(P_{2})$. Repeating such an extension process. We can finally obtain a path $P_{m-2}=v_{0}\cdots v_{m-1}$, where $v_{i}\in N_{G}(v_{i-1})\backslash V(P_{i-2})$ for $i=3,\cdots,m-1$. By this extension process, it is not difficult to check that every bipartite graph with minimum degree at least $\lfloor \frac{m}{2}\rfloor$ (or $m$) contains a path of order $m-1$ (or $2m$). $\Box$

 \begin{lem} (Hong, Liu, Lu and Ye \cite{Hong1})
Let $G$ be a graph consisting a path $P$ of order $p$ and a vertex $v \notin V(P)$. If $|N_{G}(v)\cap V(P)|\geq k\geq1$, then $G$ contains a $v$-path
of order at least $(p + k)/2 + 1$.
\end{lem}

For bipartite graphs,  the following similar result as Lemma 2.2 is obtained, which will be used to construct caterpillars and spiders in our proof.

\begin{lem}
Let $G$ be a bipartite graph consisting a path $P$ of order $p$ and a vertex $v\notin V(P)$. If $|N_{G}(v)\cap V(P)|\geq k\geq1$, then $G$ contains a $v$-path of order at least $k+(p+1)/2$.
\end{lem}

\noindent{\bf Proof.} Denote $P=v_{1} \cdots v_{p}$. Let $a=\min\{i|v_{i}\in N_{G}(v)\}$ and $b=\max\{i|v_{i}\in N_{G}(v)\}$. Since $|N_{G}(v)\cap \{v_{a}, v_{a+1}, \cdots , v_{b}\}|\geq k$ and $G$ is a bipartite graph, then $b -a + 1 =|\{v_{a}, v_{a+1}, \cdots , v_{b}\}|\geq 2k-1$. Let $P_{1} = vv_{a}Pv_{p }$ and $P_{2} = vv_{b}Pv_{1}$. Then $|V(P_{1})| + |V(P_{2})| = (p- a + 2) + (b + 1) \geq 2k+p+1$ and one of $P_{1}$, $P_{2}$ is desired. $\Box$

\begin{lem}(Hasunuma and Ono \cite{Hasunuma1})
Let $T$ be a tree of order $m$ and $S$ be a subtree of $T$. If a graph $G$ contains a subtree $S'\cong S $ such that $d_{G}(v)\geq m-1$ for any $v\in V(G)\setminus V(S')$ and for any $v\in V(S')$ with $d_{S}(\phi ^{-1}(v))<d_{T}(\phi ^{-1}(v))$, where $\phi$ is an isomorphism from $V(S)$ to $V(S')$, then $G$ contains a subtree $T'\cong T$ such that $S'\subseteq T'$.
 \end{lem}



For bipartite graphs,  we have similar result as Lemma 2.4 in the following.

\begin{lem}
Let $T=(X,Y)$ be a tree and $S$ be a subtree of $T$. If a bipartite graph $G$ contains a subtree $S'\cong S $ such that $d_{G}(v)\geq t$ for any $v\in V(G)\setminus V(S')$ and for any $v\in V(S')$ with $d_{S}(\phi ^{-1}(v))<d_{T}(\phi ^{-1}(v))$, where $\phi$ is an isomorphism from $V(S)$ to $V(S')$ and $t=\max\{|X|,|Y|\}$, then $G$ contains a subtree $T'\cong T$ such that $S'\subseteq T'$.

 \end{lem}

 \noindent{\bf Proof.} Let $G=(X_{G},Y_{G})$, $S=(X_{S},Y_{S})$, $S'=(X_{S'},Y_{S'})$, $X_{S}\subseteq X$ and $X_{S'}\subseteq X_{G}$. Without loss of generality, assume that $\phi(v)\in X_{S'}$ for any $v\in X_{S} $. Then $|X_{S}|=|X_{S'}|\leq t$ and $|Y_{S}|=|Y_{S'}|\leq t$. Let $v\in V(S')$ such that $d_{S}(\phi ^{-1}(v))<d_{T}(\phi ^{-1}(v))$. If $v\in X_{S'}$, then
\begin{equation}
\begin{aligned}
 |N_{G}(v)\setminus V(S')|
 &=d_{G}(v)-|N_{G}(v)\cap V(S')|\\
 &\geq t-|Y_{S'}|\\
 &\geq |Y|-|Y_{S}|\\
  &\geq |N_{T}(\phi^{-1}(v))\cap (Y\backslash Y_{S})|\\
  &\geq d_{T}(\phi ^{-1}(v))-d_{S}(\phi ^{-1}(v)).\\
\end{aligned}
\end{equation}
Similarly, if $v\in Y_{S'}$, then $$ |N_{G}(v)\setminus V(S')|\geq d_{T}(\phi ^{-1}(v))-d_{S}(\phi ^{-1}(v)).$$ Thus we can extend $S'$ to large subtree $U'$ of $G$ which is isomorphic to a subtree $U$ of $T$ containing $S$, by arbitrarily selecting $d_{T}(\phi ^{-1}(v))-d_{S}(\phi^{-1}(v))$ vertices in $N_{G}(v)\setminus V(S')$. Note that degree conditions for $S'$ similarly hold for $U'$. Repeating such an extension process. We can finally obtain a subtree $T'\cong T$ such that $S'\subseteq T'$. $\Box$

The set of leaves of tree $T$ is denoted by $Leaf(T)$. If $S=T-L$, where $L\subseteq Leaf(T)$, then by applying the extension process only for the neighbors of $L$ in $S$, we have the following lemma.

\begin{lem}
Let $T=(X,Y)$ be a tree with $L\subseteq Leaf(T)$ and $S=T-L$. If a bipartite graph $G$ contains a subtree $S'\cong S $ such that $d_{G}(v)\geq t$ for any $v\in V(S')$ with $d_{S}(\phi ^{-1}(v))<d_{T}(\phi ^{-1}(v))$, where $\phi$ is an isomorphism from $V(S)$ to $V(S')$ and $t=\max\{|X|,|Y|\}$, then $G$ contains a subtree $T'\cong T$ such that $S'\subseteq T'$.

\end{lem}

\section{Caterpillars}

A $caterpillar$ is a tree in which a single path (the spine) is incident to every edge. Note that paths, stars, double stars, path-stars, path-double-stars are caterpillars. Zhang \cite{Zhang} has confirmed Conjecture 2 for caterpillars when $k\leq 2$.
In this section, we will show that Conjecture 2 is true for caterpillars in 3-connected bipartite graphs.

Let $G$ be a subdivision of some simple 3-connected graph. An $ear$ of $G$ is a maximal path $P$ whose each internal vertex has degree 2 in $G$. Define $n(G) = |\{v | d_{G}(v) \geq 3\}|$. Then $G$ is 3-connected if and only if $n(G) = |V(G)|$.

\begin{thm}
For any caterpillar $T$ with bipartition $X$ and $Y$ (denote $t=\max\{|X|,|Y|\}$), every 3-connected bipartite graph $G$ with $\delta(G)\geq t+3$ contains a caterpillar $T'\cong T$ such that $G-V(T')$ is still 3-connected.
\end{thm}

\noindent{\bf Proof.} By contrary, assume that $G$ contains no caterpillar $U$ isomorphic to $T$ such that $G-V(U)$ is 3-connected. Let $T'$ be a tree in $G$ isomorphic to $T$, and let $B$ be a subdivision of some simple 3-connected graph in $G-V(T')$. The existence of $T'$ is assured by $\delta(G)\geq t+3$ and Lemma 2.5.  Since $\delta(G-V(T'))\geq 3$, $G-V(T')$ contains a subdivision of $K_{4}$. The minimum induced subgraph containing a subdivision of $K_{4}$ is a subdivision of some simple 3-connected graph. Thus $B$ also exists.  Furthermore, we assume $T'$, $B$ are such subgraphs such that $n(B)$ is as lager as possible, and $|V(B)|$ is as small as possible. We denote 
$$V_{1}=\{v\in V(B)|d_{B}(v)\geq3\}$$ and $V_{2}=V(B)\backslash V_{1} $. Let $P=u_{1}\cdots u_{r}$ be the spine of $T'$ and $H=G-V(T'\cup B)$.

If $V(H)=\emptyset $, then $\delta(B)\geq 3$ and $B$ is 3-connected, a contradiction. So $V(H)\neq\emptyset$. We will complete the proof by a series of claims.

\noindent{\bf Claim 1.} $|N_{G}(v)\cap (V(T')\backslash \{u_{1},\cdots,u_{s}\})|\leq t-\lfloor\frac{s}{2}\rfloor$ for any $v\in V(G)\backslash V(T')$ and $1\leq s\leq r$.

Since $\lfloor\frac{s}{2}\rfloor \leq |X\cap \{u_{1},\cdots,u_{s}\}|\leq \lceil\frac{s}{2}\rceil$ and $\lfloor\frac{s}{2}\rfloor \leq |Y\cap \{u_{1},\cdots,u_{s}\}|\leq \lceil\frac{s}{2}\rceil$, we have
\begin{equation}
\begin{aligned}
|N_{G}(v)\cap (V(T')\backslash \{u_{1},\cdots,u_{s}\})|&=|N_{G}(v)\cap (X\backslash  \{u_{1},\cdots,u_{s}\})|\\
&\leq |X|-\lfloor\frac{s}{2}\rfloor \\
&\leq t-\lfloor\frac{s}{2}\rfloor
\end{aligned}
\end{equation} or
 \begin{equation}
\begin{aligned}
|N_{G}(v)\cap (V(T')\backslash \{u_{1},\cdots,u_{s}\})|&=|N_{G}(v)\cap (Y\backslash \{u_{1},\cdots,u_{s}\})|\\
&\leq |Y|-\lfloor\frac{s}{2}\rfloor \\
&\leq t-\lfloor\frac{s}{2}\rfloor
\end{aligned}
\end{equation} for any $v\in V(G)\backslash V(T')$.

\noindent{\bf Claim 2.} For any $v\in V(T'\cup H)$, if $H\cup T'-\{v\}$ contains a caterpillar $T''\cong T$, then $|N_{G}(v)\cap V(B)|\leq 2$.

Suppose, to the contrary, that $v$ is such a vertex such that $|N_{G}(v)\cap V(B)|\geq 3$. If no ear of $B$ contains $N_{G}(v) \cap V (B)$, then there exist three $(v, V_1)$-paths with internal vertices not in $V(G)\backslash (V(B)\cup \{v\})$. Let $B'=G[V(B)\cup \{v\}]$. Then $B'$ is still a subdivision of some simple 3-connected graph, but $n(B') > n(B)$, a contradiction. Thus, there exists an ear $Q$ such that $N_{G}(v) \cap V (B) \subseteq V (Q)$. Let $Q = v_{1}v_{2}\cdots v_{s}$, $a=\min\{i: 1\leq i\leq s~and~ v_{i} \in N_{G}(v)\}$ and $b=\max\{i: 1\leq i\leq s ~and~ v_{i} \in N_{G}(v)\}$. Since $|N_{G}(v)\cap V(B)|\geq 3$ and $G$ is a bipartite graph, we have $s\geq 5$ and $b-a\geq 4$. Let $B'=G[V(B-\{v_{a+1},\cdots,v_{b-1}\})\cup \{v\}]$. Then $B'$ is a subdivision of some simple 3-connected graph such that $n(B')= n(B)$, but $|V(B)|-|V(B')|=b-a-2\geq 1$, also a contradiction.

\noindent{\bf Claim 3.} $|N_{G}(v)\cap V(B)|\leq 2$ for any $v\in V(H)$.

Since $T'\subseteq H\cup T'-\{h\}$ for any $h\in V(H)$, Claim 2 implies Claim 3.

\noindent{\bf Claim 4.} $|N_{G}(u_{i})\cap V(B)|\leq 2$ for each $i\in\{1,\cdots,r\}$.

Suppose that $|N_{G}(u_{i})\cap V(B)|\geq 3$ for some $i\in\{1,\cdots,r\}$. Let $l=\min\{i:1\leq i\leq r ~and~ |N_{G}(u_{i})\cap V(B)|\geq 3 \}.$ Then \begin{equation}
\begin{aligned}|N_{G}(v)\cap(V(G)\setminus (V(B)\cup\{u_{l}\}))|= d_{G}(v)-|N_{G}(v)\cap V(B)|-|N_{G}(v)\cap \{u_{l}\}|\geq t\\\end{aligned}
\end{equation} for any $v\in V(H)\cup\{u_{1},\cdots,u_{l-1}\}$, and $|N_{G}(u_{i})\cap V(H)|= d_{G}(u_{i})-|N_{G}(u_{i})\cap V(B)|-|N_{G}(u_{i})\cap V(T)|\geq 1$ for each $i\in\{1,\cdots,l-1\}$.


Suppose $N_{G}(u_{i})\cap(V(H)\cup\{u_{1},\cdots,u_{l-1}\})=\emptyset$ for each $i\in \{l+1, \cdots, r\}$, then\begin{equation}
\begin{aligned}|N_{G}(h)\backslash V(B\cup T')|&= d_{G}(h)-|N_{G}(h)\cap V(B)|\\
&-|N_{G}(h)\cap (V(T')\backslash \{u_{l+1},\cdots,u_{r}\})|\\
&\geq t+1-(t-\lfloor\frac{r-l}{2}\rfloor)\\
&\geq 1+\lfloor\frac{r-l}{2}\rfloor,\end{aligned}
\end{equation} for all $h\in V(H)$. Thus, $\delta(H)\geq1+\lfloor\frac{r-l}{2}\rfloor$ and $H$ contains a path of order at least $(r-l+1)$ by Lemma 2.1. Since $|N_{G}(u_{l-1})\cap V(H)|\geq 1$, there exits a path $Q'=h_{1}\cdots h_{r-l+1}$ in $H$ such that $h_{1}\in N_{G}(u_{l-1})$. Let $S'=Q'\cup R'\cup \{u_{l-1}h_{1}\}$ be a tree, where $R'=T'[\{u_{1},\cdots,u_{l-1}\}\cup (N _{T'}(\{u_{1},\cdots,u_{l-1}\})\setminus \{u_{l}\})] $. By Lemma 2.6 and inequality (4), there exists a caterpillar $T''\cong T$ in $G-V(B)\cup\{u_{l}\}$, a contradiction to Claim 2. Thus, $N_{G}(u_{i})\cap(V(H)\cup\{u_{1},\cdots,u_{l-1}\})\neq\emptyset$ for some $i\in \{l+1, \cdots, r\}$. Assume $j=\min\{i: l+1\leq i\leq r~ and~ |N_{G}(u_{i})\cap(V(H)\cup\{u_{1},\cdots,u_{l-1}\})| \geq1 \}.$ Let us consider the following two cases.


\noindent{\bf Case 1.} $j=l+1$.

For any $v\in V(H)\cup\{u_{1},\cdots,u_{l-1}\} $,
  \begin{equation}
\begin{aligned}
  |N_{G}(v)\cap(V(H)\cup\{u_{1},\cdots,u_{l-1}\})|&\geq d_{G}(v)-|N_{G}(v)\cap V(B)|\\
  &-|N_{G}(v)\cap (V(T')\setminus\{u_{1},\cdots,u_{l-1}\})|\\
  &\geq (t+3)-2-(t-\lfloor \frac{l-1}{2}\rfloor)\\
  &\geq \lfloor \frac{l-1}{2}\rfloor+1.\end{aligned}
\end{equation} 
Thus, by Lemma 2.1, $G[ V(H)\cup\{u_{1},\cdots,u_{l-1}\}]$ contains a path $Q_{1}=w_{1}\cdots w_{l}$ such that $w_{l}\in N_{G}(u_{l+1})$. Let $S_{1}=Q_{1}\cup R_{1}\cup \{w_{l}u_{l+1}\}$ be a tree, where $R_{1}=T'[ \{u_{l+1},\cdots,u_{r}\}\cup (N _{T'}(\{u_{l+1},\cdots,u_{r}\})\setminus \{u_{l}\})] $. By Lemma 2.6 and inequality (4), there exists a caterpillar $T_{1}\cong T$ in $G-V(B)\cup\{u_{l}\}$ such that $S_{1}\subseteq T_{1}$, a contradiction to Claim 2.

\noindent{\bf Case 2.} $j>l+1$.

By $j>l+1$, we have $N_{G}(v)\cap \{u_{l+1},\cdots,u_{j-1}\}=\emptyset$ for each $v\in V(H)\cup\{u_{1},\cdots,u_{l-1}\}$ and $j\geq3$. We will find a $u_{j}$-path of order at least $j$ in $G[\{u_{1},\cdots,u_{l-1},u_{j}\}\cup V(H)]$ according to the following two subcases.

\noindent{\bf Subcase 2.1} $N_{G}(u_{j})\cap V(H)=\emptyset$.

By $N_{G}(u_{j})\cap V(H)=\emptyset$, we have $N_{G}(u_{j})\cap \{u_{1},\cdots,u_{l-1}\}\neq \emptyset$ and
\begin{equation}
\begin{aligned}
|N_{G}(h)\backslash V(T'\cup B)|&\geq d_{G}(h)-|N_{G}(h)\cap V(B)|\\
&-|N_{G}(h)\cap(V(T')\backslash\{u_{l+1},\cdots,u_{j}\})|\\
&\geq (t+3)-2-(t-\lfloor\frac{j-l}{2}\rfloor)\\
&\geq1+\lfloor\frac{j-l}{2}\rfloor,
\end{aligned}
\end{equation} for each $h\in V(H)$. Thus $\delta (H)\geq 1+\lfloor\frac{j-l}{2}\rfloor$. Let $a=\min\{i: 1\leq i\leq l-1 ~and~ u_{i}\in N_{G}(u_{j})\}$ and $b=\max\{i:1\leq i\leq l-1~and~ u_{i}\in N_{G}(u_{j})\}.$ Then $b\geq a $. Let $P_{1} = u_{j}u_{a}u_{a+1}\cdots u_{l-1}$ and $P_{2} = u_{j}u_{b}u_{b-1}\cdots u_{1}$. Then $|V(P_{1})|<j$ and $|V(P_{2})|<j$.

Suppose $|V(P_{1})|\geq |V(P_{2})|$, then $|V(P_{1})|=l-a+1\geq \lfloor\frac{l-1}{2}\rfloor+2$. Since $|N_{G}(u_{l-1})\cap V(H)|\geq 1$, there exits a path $P_{H}=h_{1}\cdots h_{p}$ in $H$ such that $h_{p}\in N_{G}(u_{l-1})$ and $p\geq j-l+1$. If $ l-a+1+p\geq j$, then there is a $u_{j}$-path of order at least $j$ in $G[\{u_{a},\cdots,u_{l-1},u_{j}\}\cup V(H)]$. Otherwise, \begin{equation}
\begin{aligned}
|N_{G}(h_{p})\cap\{u_{1},\cdots,u_{a-1}\})|&\geq d_{G}(h_{p})-|N_{G}(h_{p})\cap V(H)|-|N_{G}(h_{p})\cap V(B)|\\
&-|N_{G}(h_{p})\cap (V(T')\backslash\{u_{1},\cdots, u_{a-1},u_{l+1},\cdots, u_{j}\})|\\
&\geq (t+3)-\lceil\frac{p-1}{2}\rceil-2-(t-\lfloor\frac{a-1}{2}\rfloor-\lfloor\frac{j-l}{2}\rfloor)\\
&=\lfloor\frac{a-1}{2}\rfloor+\lfloor\frac{j-l}{2}\rfloor+1-\lceil\frac{p-1}{2}\rceil\\
&\geq \lfloor\frac{j-l+a-1}{2}\rfloor-\lceil\frac{p-1}{2}\rceil\\
&>0,\\
\end{aligned}
\end{equation}
and there is a $u_{j}$-path of order at least \begin{equation}
\begin{aligned}
(\lfloor\frac{a-1}{2}\rfloor+\lfloor\frac{j-l}{2}\rfloor+1&-\lceil\frac{p-1}{2}\rceil)+\frac{a}{2}+l-a+p\\
&\geq \lfloor\frac{j-l}{2}\rfloor+l+\lceil\frac{p}{2}\rceil\\
&\geq \lfloor\frac{j-l}{2}\rfloor+l+\lceil\frac{j-l+1}{2}\rceil\\
&\geq j,\\
\end{aligned}
\end{equation} in $G[\{u_{1},\cdots,u_{l-1},u_{j}\}\cup V(H)] $ by Lemma 2.3.

Similarly, for $|V(P_{2})|\geq |V(P_{1})|$, we also can find  a $u_{j}$-path of order at least $j$ in $G[\{u_{1},\cdots,\\ u_{l-1},u_{j}\}\cup V(H)]$.

\noindent{\bf Subcase 2.2} $N_{G}(u_{j})\cap V(H)\neq\emptyset$.

For any $h\in V(H)$, we have
\begin{equation}
\begin{aligned}
|N_{G}(h)\backslash V(T'\cup B)|&\geq d_{G}(h)-|N_{G}(h)\cap V(B)|\\
&-|N_{G}(h)\cap(V(T')\backslash\{u_{l+1},\cdots,u_{j-1}\}|\\
&\geq (t+3)-2-(t-\lfloor\frac{j-l-1}{2}\rfloor)\\
&\geq1+\lfloor\frac{j-l-1}{2}\rfloor.
\end{aligned}
\end{equation}
Thus $\delta (H)\geq 1+\lfloor\frac{j-l-1}{2}\rfloor$ and $H$ contains a path of order at least $(j-l)$ by Lemma 2.1. Let $P'_{H}=z_{1}\cdots z_{q}$ be a longest path in $H$ such that $z_{q}\in N_{G}(u_{j})$. Then  $q\geq j-l$ and 
\begin{equation}
\begin{aligned}
|N_{G}(z_{1})\cap\{u_{1},\cdots,u_{l-1}\})|&\geq d_{G}(z_{1})-|N_{G}(z_{1})\cap V(H)|-|N_{G}(z_{1})\cap V(B)|\\
&-|N_{G}(z_{1})\cap \{u_{l}\}|-|N_{G}(z_{1})\cap (V(T')\backslash\{u_{1},\cdots, u_{j-1}\})|\\
&\geq (t+3)-\lceil\frac{q-1}{2}\rceil-2-1-(t-\lfloor\frac{j-1}{2}\rfloor)\\
&\geq\lfloor\frac{j-1}{2}\rfloor-\lceil\frac{q-1}{2}\rceil.\\
\end{aligned}
\end{equation}
If $N_{G}(z_{1})\cap\{u_{1},\cdots,u_{l-1}\}=\emptyset$, then $q\geq j-1$ and there exits a $u_{j}$-path of order at least $j$ in $G[\{u_{j}\}\cup V(H)]$. Otherwise, $N_{G}(z_{1})\cap\{u_{1},\cdots,u_{l-1}\})\neq \emptyset$ and there exists a $u_{j}$-path of order at least $j$ in $G[\{u_{1},\cdots,u_{l-1},u_{j}\}\cup V(H)]$ by Lemma 2.3 and \begin{equation}
\begin{aligned}\lfloor\frac{j-1}{2}\rfloor-\lceil\frac{s-1}{2}\rceil+\frac{l}{2}+q+1&\geq \lceil\frac{j}{2}\rceil+\frac{l}{2}+\lfloor\frac{q}{2}\rfloor\\
&\geq \lceil\frac{j}{2}\rceil+\frac{l}{2}+\lfloor\frac{j-l}{2}\rfloor\\
&\geq j-\frac{1}{2}.
\end{aligned}
\end{equation}

Thus, there exists a $u_{j}$-path of order at least $j$ in $G[\{u_{1},\cdots,u_{l-1},u_{j}\}\cup V(H)]$. Let $S_{2}=Q_{2}\cup R_{2}\cup \{\omega_{j-1}u_{j}\}$ be a tree, where $Q_{2}=w_{1}\cdots w_{j-1}$ is a path in $G[\{u_{1},\cdots,u_{l-1},u_{j}\}\cup V(H)]$, $w_{j-1}\in N_{G}(u_{j})$ and $R_{2}=T'[ \{u_{j},\cdots,u_{r}\}\cup (N _{T'}\{u_{j},\cdots,u_{r}\}\setminus \{u_{j-1}\})] $. By Lemma 2.6 and inequality (4), there exists a caterpillar $T_{2}\cong T$ in $G-V(B)\cup\{u_{l}\}$ such that $S_{2}\subseteq T_{2}$, a contradiction to Claim 2. Thus Claim 4 holds.

\noindent{\bf Claim 5.} $|N_{G}(v)\cap V(B)|\leq 2$ for each $v\in V(T')\backslash V(P)$.

 Suppose, that there exists $w\in V(T')\backslash V(P)$ such that $|N_{G}(w)\cap V(B)|\geq 3$ and $w u_{i}\in E(T)$ for some $i$. By Claim 4 and $$|N_{G}(u_{i})\setminus V(T'\cup B)|= d_{G}(u_{i})-|N_{G}(u_{i})\cap V(T')|-|N_{G}(u_{i})\cap V( B)|\geq 1,$$ then there exists a vertex $w'\in V(H)$ such that $u_{i}w'\in E(G)$. We can find the caterpillar $T''\cong T$ from $T'$ by deleting $w$ and adding $u_{i}w'$, a contradiction to Claim 2.

\noindent{\bf Claim 6.} $V_{2}\neq \emptyset$.

Assume, to the contrary, that $V_{2}=\emptyset$. Then $\delta (B)\geq 3$ and $B$ is 3-connected. Since $G$ is a 3-connected graph, there exist three $(v,B)$-paths $P_{1}$, $P_{2}$, $P_{3}$  for any $v\in G-V(B)$. We choose $v$ and  $P_{1}$, $P_{2}$, $P_{3}$ such that $\sum^{3}_{i=1}|V(P_{i})|$ is as small as possible. Let $B'=G[V(B)\cup V(P_{1}\cup P_{2} \cup P_{3})]$. Then $B'$ is still a subdivision of some simple 3-connected graph by $\delta (B)\geq 3$. If $\delta (G-V(B'))\geq t$, then there exists a caterpillar $T''\cong T$ in $G-V(B')$, but $n(B')>n(B)$, which is a contradiction. If $\delta (G-V(B'))\leq t-1$, then there exists a vertex $w\in V(G) \backslash V(B')$ such that $N_{G}(w)\cap (V(G)\setminus V(B'))\leq t-1$. Thus $|N_{G}(w)\cap V(B')|=d_{G}(w)-|N_{G}(w)\cap (V(G)\setminus V(B')|\geq 4$. Let $v_{i}\in N_{G}(w)\cap V(B')$ such that $v_{i}\notin \{v_{1}, \cdots, v_{i-1}\}$ and the distance from $v_{i}$ to $B$ in $B'-\{v_{1}, \cdots, v_{i-1}\}$ is as small as possible for $i=1,2,3$. Then, we can obtain three $(w, V (B))$-paths $P'_{1}$, $P'_{2}$, $P'_{3}$ in $G[V (B')\cup\{w\}]$ using edges $wv_{1}, wv_{2}, wv_{3}$ such that $\sum^{3}_{i=1}|V(P_{i})|> \sum^{3}_{i=1}|V(P'_{i})|$, a contradiction. Thus Claim 6 holds.

Since $V_{2}\neq \emptyset$, there exists a shortest $(V_{2},B)$-path $Q$ which is edge-disjoint with $B$. Let $B'=G[V(B)\cup V(Q)]$. Then $B'$ is still a subdivision of some simple 3-connected graph and $n(B')>n(B)$. Suppose $\delta (G-V(B'))\geq t$, then $G-V(B')$ contains a caterpillar isomorphism to $T$ by Lemma 2.5, but $n(B')>n(B)$, which is a contradiction. Thus $\delta (G-V(B'))\leq t-1$, and then there exists a vertex $v\in V(G)\backslash V(B')$ such that $|N_{G}(v)\cap V(B')|=d_{G}(v)-|N_{G}(v)\cap (V(G)\setminus V(B')|\geq 4$. By Claims 3-5, we have $|N_{G}(v)\cap V(B)|\leq 2$ and $|N_{G}(v)\cap (V(Q)\setminus End(Q))|\geq 2$. Since $G$ is a bipartite graph and $Q$ is a shortest $(V_{2},B)$-path, we further have  $|N_{G}(v)\cap V(B)|=2$ and $|N_{G}(v)\cap (V(Q)\setminus End(Q))|=2$. Let $End(Q)=\{u,u'\}$, where $u\in V_{2}$ and $u'\in V(B)$. Assume that $v_{1}$ is a neighbor of $v$ nearest to $u$ in $V(Q)\setminus End(Q)$ and $w_1$ is a neighbor of $v$  in $V(B)$. Then $Q'=uQv_{1}vw_{1}$ is a shorter $(V_{2},B)$-path than $Q$, which is  a contradiction to the choice of $Q$. The proof is thus complete. $\Box$

\section{Spiders}

In this section, we will confirm Conjecture 2 for spiders when $k\leq 3$. A $spider$ is a tree with at most one vertex of degree more than 2, called the center of the spider (if no vertex of degree more than two, then any vertex can be the center). A $leg$ of a spider is a path from the center to a vertex of degree 1. Then each leg of a spider with center $u$ is a $u$-path. Usually, all the legs of a spider are non-trivial $u$-paths.

\begin{thm}
For every spider $T$ with bipartition $X$ and $Y$, every 2-connected bipartite graph $G$ with $\delta(G)\geq t+2$ contains a spider $T' \cong T$ such that $G-V(T')$ is still 2-connected, where $t=\max \{|X|, |Y|\}$.

\end{thm}

\noindent{\bf Proof.} By Lemma 2.5, there exists a subgraph in $G$ isomorphic to $T$. Then we assume $T'\cong T$ is such a subgraph of $G$ such that the maximum block $B$ of $G-V(T')$ has order as lager as possible. Assume $T'$ has center $v$ and has $d$ legs, namely $P_{1}, P_{2},\cdots , P_{d}$ whose length are $p_{1}, p_{2},\cdots , p_{d}$,  respectively. Also, assume that $p_{i} \geq 2$ for $i = 1, \cdots ,r $ and $p_{i} = 1$ for $i = r + 1, \cdots , d$.

 If $G-V(T')$ is 2-connected, then we are done. So we assume $H=G-V(B\cup T')$ is not empty.
By the assumption of $B$, then $|N_{G}(h)\cap V(B)|\leq 1$ for each $h\in V(H)$. Otherwise, $G-V(T')$ has a block contains $V(B)\cup \{h\}$, a contradiction to the choice of $B$.
We denote
$$V_{1}=\{u\in V(T')| |N_{G}(u)\cap V(B)|\geq 2\}$$
and $V_{2}=N_{G}(H)\cap V(T')$. We have the following claims.

\textbf{Claim 1.} $V_{1}\cup V_{2}=V(T')$, $V_{1}\neq \emptyset $ and $V_{2}\neq \emptyset $.

Suppose that there exists $u\in V(T')\setminus (V_{1}\cup V_{2})$, then $|N_{G}(u)\cap V(B)|\leq 1$ and $N_{G}(u)\cap V(H)=\emptyset$ by the definition of $V_{1}$ and $V_{2}$. Thus $d_{G}(u)=|N_{G}(u)\cap V(B)|+|N_{G}(u)\cap V(H)|+|N_{G}(u)\cap V(T')| \leq 1+t$, a contradiction to the assumption of the minimum degree.

If $V_{1}=\emptyset$, then $|N_{G}(u)\cap V(B)|\leq 1$ for each $u\in V(H\cup T')$. Since $G$ is 2-connected, there is a shortest path $Q$ with both ends in $B$ and inner vertices not in $B$ such that $|V(Q)|\geq 4$. Let $Q=v_{1} \cdots v_{s}$ and $v_{1},v_{s}\in V(B)$. Then $N_{G}(v_{2})\cap V(B)=\{v_{1}\}$ and $N_{G}(v_{2})\cap (V(B\cup Q))=\{v_{1},v_{3}\}$. Thus $|N_{G}(v_{2})\setminus V(B\cup Q)|\geq t $ and $V(G)\setminus V(B\cup Q)\neq \emptyset$. By the assumption of $Q$, $|N_{G}(w)\cap V(Q)|\leq 2$ for each $w\notin V(B\cup Q)$. Thus $|N_{G}(w)\cap V(B\cup Q)|\leq 2$ and $\delta (G-(V(B\cup Q)))\geq t$. Hence $G-V(B\cup Q)$ contains a spider $T''\cong T$ by Lemma 2.5. However, $G-V(T'')$ has a block contains $V(B)\cup V(Q)$, which is a contradiction.

If $V_{2}= \emptyset $, then $N_{G}(H)\cap V(T')= \emptyset$. Since $G$ is 2-connected, we have $N_{G}(H)\cap V(B)\geq 2$. Thus there exists a path $Q'$ with both ends in $B$ and inner vertices in $H$. Let $B'=G[V(B\cup Q')]$. Then $B'$ is a block in $G-V(T')$ with $|V(B')|>|V(B)|$, which is a contradiction to the choice of $B$. Thus Claim 1 holds.

We consider two cases in the following proof.

\textbf{Case 1.} $v\in V_{1}$.

For any $h \in V(H)$, denote by $I_{h} = \{i| i \leq r ~and~ N_{G}(h) \cap End(P_{i} -v) \neq \emptyset\}$. Pick $u \in V(H)$ so that $|I_{u}|$ is as large as possible. Without loss of generality, we assume that $I_{u} = \{1,\cdots, s\}$ for some $s\leq r$ and $|V(P_{s+1})|\geq\cdots\geq |V(P_{r})|$. We will construct a spider $T' _{1}$ in $H$ with center $u$ and with $r$ legs ($P'_{1}, P'_{2}, \cdots, P'_{r}$) such that $|V(P'_{i})|=p_{i}$ for any $i\in\{1,\cdots, r\}$. If $i\in\{1,\cdots, s\}$, then $P'_{i}$ is not difficult to construct by the choice of $u$. Thus, we assume  $i\in\{s+1,\cdots, r\}$.

Let $Q_{s+1}'$ be the longest $u$-path in $H-\bigcup_{i=1}^{s}V(P_{i}'-u)$ and $z'_{s+1}$ be the end of $Q'_{s+1}$. If $|V(Q_{s+1}')|\geq p_{s+1}$, then $Q'_{s+1}$ contains a $u$-path of order at least $ p_{s+1}$, for otherwise $|V(Q_{s+1}')|< p_{s+1}$. If $s+1\in I_{z'_{s+1}}$, then $G[V(Q_{s+1}')\cup (V(P_{s+1})-\{v\})]$ contains a $u$-path of order at least $ p_{s+1}$. If $s+1\notin I_{z'_{s+1}}$ then\begin{equation}
\begin{aligned} |N_{G}(z'_{s+1})\cap (V(P_{s+1})-\{v\})|
&=d_{G}(z'_{s+1})-|N_{G}(z'_{s+1})\cap V(B)|-|N_{G}(z'_{s+1})\cap V(H)|\\
&-|N_{G}(z'_{s+1})\cap (V(T')\setminus (V(P_{s+1})-\{v\}))|\\
&\geq (t+2)-1-\lceil\frac{|V(Q_{s+1}')|-1}{2}\rceil-(t-\lfloor\frac{|V(P_{s+1})|-1}{2}\rfloor)\\
&\geq 1-\lceil\frac{|V(Q_{s+1}')|}{2}\rceil+\lfloor\frac{|V(P_{s+1})|-1}{2}\rfloor\\
&\geq\lceil\frac{|V(P_{s+1})|}{2}\rceil-\lceil\frac{|V(Q_{s+1}')|}{2}\rceil.\\
\end{aligned}
\end{equation}
If $N_{G}(z'_{s+1})\cap (V(P_{s+1})-\{v\})=\emptyset$, then $|V(Q_{s+1}')|\geq p_{s+1}$. Otherwise, there is a $u$-path of order at least $$(\lceil\frac{|V(P_{s+1})|}{2}\rceil-\lceil\frac{|V(Q_{s+1}')|}{2}\rceil+\frac{|V(P_{s+1})|-1+1}{2})+|V(Q_{s+1}')|\geq p_{s+1}$$
in $G[V(Q_{s+1}')\cup (V(P_{s+1})-\{v\})]$ by Lemma 2.3.

Repeating such a process. We can finally obtain a spider $T '_{1}$ with $r$ legs whose lengths are $p_{1}, \cdots , p_{r}$, respectively. Since $|N_{G}(u)\setminus (V(B)\cup \{v\})|\geq t$, then $G[V(G)\backslash (V(B)\cup \{v\})]$ contains a spider $T'_{2}\cong T$ such that $T_{1}'\subseteq T'_{2}$ by Lemma 2.6 and
$V(B)\cup \{v\}$ is contained in a block, a contradiction to the choice of $T'$.

\textbf{Case 2.} $v\notin V_{1}$.

By $v\notin V_{1}$, we may find an edge $w_{1}w_{2}$ in some leg $P_{i}$ ($i\in \{1,\cdots, d\}$) such that $w_{1} \in V_{2}$, $w_{2} \in V_{1}$, $ w_{1} \in V(vP_{i}w_{2})$ and subject to this, $vP_{i}w_{2}$ is as long as possible. Let $w$ be the end of $P_{i}$ and $Q_{i} = w_{2}P_{i}w- \{w_{2}\}$. Let $w' \in N_{G}(w_{1}) \cap V(H)$ and $Q'_{i}$ be the longest $w'$-path in $H$. Assume $w''$ is the end of $Q'_{i}$. Then\begin{equation}
\begin{aligned} |N_{G}(w'')\cap V(Q_{i})|
&=d_{G}(w'')-|N_{G}(w'')\cap V(B)|-|N_{G}(w'')\cap V(H)|\\
&-|N_{G}(w'')\cap (V(T')\setminus V(Q_{i}))|\\
&\geq (t+2)-1-\lceil\frac{|V(Q'_{i})|-1}{2}\rceil-(t-\lfloor\frac{|V(Q_{i})|}{2}\rfloor)\\
&\geq 1+\lfloor\frac{|V(Q_{i})|}{2}\rfloor-\lceil\frac{|V(Q'_{i})|}{2}\rceil.\\
\end{aligned}
\end{equation}
If $N_{G}(w'')\cap V(Q_{i})=\emptyset$, then $|V(Q'_{i})|\geq |V(Q_{i})|+1$. Otherwise, there is a $w'$-path of order at least $$(1+\lceil\frac{|V(Q_{i})|}{2}\rceil-\lceil\frac{|V(Q'_{i})|}{2}\rceil+\frac{|V(Q_{i})|+1}{2})+|V(Q'_{i})|\geq |V(Q_{i})|+1$$
in $G[V(Q_{i}'\cup Q_{i})]$ by Lemma 2.3. Then $V(H)\cup V(T'-\{w_{2}\})$ contains a spider $T'_{3} \cong T'$ such that $G[V(B)\cup \{w_{2}\}]$ is contained in a block of $G-V(T')$, which is a contradiction. The proof is finished. $\Box$

The proof of Theorem 4.2 is similar to that of Theorem 3.1.  In fact, we only need to  replace the method of constructing caterpillars in Theorem 3.1 with the method of constructing spiders in Theorem 4.1. So we omit the detail proof here.

\begin{thm}
For every spider $T$ with bipartition $X$ and $Y$, every 3-connected bipartite graph $G$ with $\delta(G)\geq t+3$ contains a spider $T' \cong T$ such that $G-V(T')$ is still 3-connected, where $t=\max \{|X|, |Y|\}$.

\end{thm}

%


\begin{thebibliography}{5}

\bibitem{Bondy} J. A. Bondy, U. S. R. Murty, Graph Theory, Graduate Texts in Mathematics 244, Springer, Berlin, 2008.

\bibitem{Chartrand} G. Chartrand, A. Kaigars, D.R. Lick, Critically n-connected graphs, Proc. Amer. Math. Soc. 32 (1972) 63-68.

\bibitem{Diwan} A. A. Diwan, N. P. Tholiya, Non-separating trees in connected graphs, Discrete Math. 309 (16) (2009) 5235-5237.

\bibitem{Fujita} S. Fujita, K. Kawarabayashi, Connectivity keeping edges in graphs with large minimum degree, J. Combin. Theory Ser. B 98 (2008) 805-811.

\bibitem{Hasunuma1} T. Hasunuma, K. Ono, Connectivity keeping trees in 2-connected graphs, J. Graph Theory 94 (2020) 20-29.

\bibitem{Hasunuma2} T. Hasunuma, Connectivity Keeping Trees in 2-Connected Graphs with Girth Conditions, Algorithmica 83 (2021) 2697-2718.

\bibitem{Hong1} Y. Hong, Q. Liu, C. Lu, Q. Ye, Connectivity keeping caterpillars and spiders in 2-connected graphs, Discrete Math. 344 (3) (2021) 112236.

\bibitem{Hong} Y. Hong, Q. Liu, Mader's conjecture for graphs with small connectivity, J. Graph Theory (2022) 1-10, https://doi.org/10.1002/jgt.22831.


\bibitem{Lu} C. Lu, P. Zhang, Connectivity keeping trees in 2-connected graphs, Discrete Math. 343 (2)(2020) 1-4.

\bibitem{Luo} L. Luo, Y. Tian, L. Wu, Connectivity keeping paths in $k$-connected bipartite graphs, Discrete Math. 345 (4) (2022) 112788.

\bibitem{Mader1} W. Mader, Connectivity keeping paths in $k$-connected graphs, J. Graph Theory 65 (2010) 61-69.

\bibitem{Mader2} W. Mader, Connectivity keeping trees in $k$-connected graphs, J. Graph Theory 69 (2012) 324-329.

\bibitem{Tian1} Y. Tian, J. Meng, H. Lai, L. Xu, Connectivity keeping stars or double-stars in 2-connected graphs, Discrete Math. 341 (4) (2018) 1120-1124.

\bibitem{Tian2} Y. Tian, H. Lai, L. Xu, J. Meng, Nonseparating trees in 2-connected graphs and oriented trees in strongly connected digraphs, Discrete Math. 342 (2) (2019) 344-351.

\bibitem{Zhang} P. Zhang, Research on the Existence of Subgraphs Related to Connectivity (Ph.D. dissertation), East China Normal University, 2021.

\end{thebibliography}
\end{document}